# On the Conditional Distribution of a Multivariate Normal given a Transformation - the Linear Case

BY RAJESHWARI MAJUMDAR* AND SUMAN MAJUMDAR

We show that the orthogonal projection operator onto the range of the adjoint of a linear operator $T$ can be represented as $UT$, where $U$ is an invertible linear operator. Given a Normal random vector $Y$ and a linear operator $T$, we use this representation to obtain a linear operator $\widehat{T}$ such that $\widehat{T}Y$ is independent of $TY$ and $Y - \widehat{T}Y$ is an affine function of $TY$. We then use this decomposition to prove that the conditional distribution of a Normal random vector $Y$ given $\mathcal{T}Y$, where $\mathcal{T}$ is a linear transformation, is again a multivariate Normal distribution. This result is equivalent to the well-known result that given a $k$-dimensional component of a $n$-dimensional Normal random vector, where $k < n$, the conditional distribution of the remaining $(n-k)$-dimensional component is a $(n-k)$-dimensional multivariate Normal distribution, and sets the stage for approximating the conditional distribution of $Y$ given $g(Y)$, where $g$ is a continuously differentiable vector field.

**1. Introduction.** What can we ascertain about the conditional distribution of a multivariate Normal random vector $Y \in \Re^n$ given $g(Y)$, where $g : \Re^n \mapsto \Re^m$ is a measurable function? Clearly, given a particular functional form of $g$, the problem is a very specific one, and depending on the functional form, may or may not have a closed form solution. Our objective is to derive an approximation to the conditional distribution in question based on some regularity properties of $g$. Specifically, in this paper we find the conditional distribution when $g$ is a linear transformation, and in a companion paper use that to derive the desired approximation when $g$ is a continuously differentiable ($C^1$) vector field by exploiting the *local linearity* of $g$.

Before proceeding further, we present a brief review of what is known about the conditional distribution when $g$ is a linear transformation, to be denoted by $\mathcal{T}$ in what follows. Casella and Berger (2002) define the bivariate Normal distribution by specifying the joint density in terms of the five parameters of the distribution - the means $\mu_1$ and $\mu_2$, the variances $\sigma_1^2$ and $\sigma_2^2$, and the correlation $\rho$ [Definition 4.5.10]. They calculate the marginal density of $Y_1$ and note that (using the joint and marginal densities) it is easy to verify that the conditional distribution of $Y_2$ given $Y_1 = y_1$ is Normal with

$$\text{mean} = \mu_2 + \rho\frac{\sigma_2}{\sigma_1}(y_1 - \mu_1) \text{ and variance} = \sigma_2^2\left(1 - \rho^2\right).$$

Anderson (1984, Section 2.5.1) and Flury (1997, Theorem 3.3.1) generalize this result to

---

the multivariate Normal distribution. While Anderson (1984) deals with the Lebesgue density of the multivariate Normal distribution (which exists only when the covariance matrix is of full rank), Flury (1997) avoids dealing with the density by defining the multivariate Normal distribution in terms of linear functionals, but requires the covariance matrix of the conditioning component to be of full rank.

Though their approaches to defining the multivariate Normal distribution differ, both Muirhead (1982, Theorem 1.2.11) and Eaton (1983, Proposition 3.13) obtain, as described below, the conditional distribution without any restriction on the rank of the covariance matrix. Let $Y_{n\times 1}$ have the multivariate Normal distribution with mean vector $\mu_{n\times 1}$ and covariance matrix $\Sigma_{n\times n}$. Partition $Y$, $\mu$, and $\Sigma$ as

$$Y = \begin{pmatrix} Y_1 \\ Y_2 \end{pmatrix}, \quad \mu = \begin{pmatrix} \mu_1 \\ \mu_2 \end{pmatrix}, \quad \Sigma = \begin{bmatrix} \Sigma_{11} & \Sigma_{12} \\ \Sigma_{21} & \Sigma_{22} \end{bmatrix}$$

where $Y_1$ and $\mu_1$ are $k \times 1$ and $\Sigma_{11}$ is $k \times k$. Then $Y_2 - \Sigma_{21}\Sigma_{11}^{-}Y_1$ and $Y_1$ are jointly Normal and uncorrelated, hence independent, where $\Sigma_{11}^{-}$ is a generalized inverse of $\Sigma_{11}$. Consequently, the conditional distribution of $Y_2 - \Sigma_{21}\Sigma_{11}^{-}Y_1$ given $Y_1$ equals the unconditional distribution of $Y_2 - \Sigma_{21}\Sigma_{11}^{-}Y_1$, which is multivariate Normal in $(n-k)$ dimensions with mean $\mu_2 - \Sigma_{21}\Sigma_{11}^{-}\mu_1$ and covariance $\Sigma_{22} - \Sigma_{21}\Sigma_{11}^{-}\Sigma_{12}$. Thus, the conditional distribution of $Y_2$ given $Y_1$ is multivariate Normal in $(n-k)$ dimensions with mean $\mu_2 + \Sigma_{21}\Sigma_{11}^{-}(Y_1 - \mu_1)$ and covariance $\Sigma_{22} - \Sigma_{21}\Sigma_{11}^{-}\Sigma_{12}$.

Given two topological vector spaces $V$ and $W$, let $\mathcal{L}(V, W)$ denote the linear space of continuous linear transformations from $V$ to $W$, and let $\mathcal{L}(V, V) = \mathcal{L}(V)$. Now note that, if $\mathcal{T} \in \mathcal{L}(\Re^n, \Re^m)$, the joint distribution of $(\mathcal{T}Y, Y)' \in \Re^{m+n}$ is multivariate Normal with

$$\text{mean} = \begin{pmatrix} A\mu \\ \mu \end{pmatrix} \text{ and covariance} = \begin{bmatrix} A\Sigma A' & A\Sigma \\ \Sigma A' & \Sigma \end{bmatrix}$$

[Muirhead (1982, Theorem 1.2.6)], where the $m \times n$ matrix $A$ represents the transformation $\mathcal{T}$ with respect to the standard orthonormal bases of $\Re^n$ and $\Re^m$. By the results of Muirhead (1982) and Eaton (1983) cited in the previous paragraph, we obtain that the conditional distribution of $Y$ given $\mathcal{T}Y$ is multivariate Normal in $\Re^n$ with mean $\mu + \Sigma A'(A\Sigma A')^{-}A(Y - \mu)$ and covariance $\Sigma - \Sigma A'(A\Sigma A')^{-}A\Sigma$.

Unfortunately, this derivation of the conditional distribution of $Y$ given $\mathcal{T}Y$, because of its dependence on manipulative matrix algebra, is not of much help when it comes to approximating the conditional distribution of $Y$ given $g(Y)$ for a $g$ in $C^1$ by exploiting the local linearity of $g$. In what follows, we present an alternative derivation of the conditional distribution of $Y$ given $\mathcal{T}Y$. Given $T \in \mathcal{L}(\Re^n)$, we find in Theorem 3 $\widehat{T} \in \mathcal{L}(\Re^n)$, depending on $T$ and the covariance of $Y$, such that $\widehat{T}Y$ is independent of $TY$ and $Y - \widehat{T}Y$ is an affine function of $TY$. In Theorem 4, given $\mathcal{T} \in \mathcal{L}(\Re^n, \Re^m)$, we find $T \in \mathcal{L}(\Re^n)$ such that the conditional distribution of $Y$ given $TY$ equals that given $\mathcal{T}Y$, and use the decomposition obtained in Theorem 3 to obtain the conditional distribution of $Y$ given $TY$, hence that given $\mathcal{T}Y$. Our derivation facilitates the approximation of the conditional distribution when the transformation $g$ is nonlinear but continuously differentiable; see Remark 4 for a brief outline.

We define a multivariate Normal distribution in terms of an affine transformation of the random vector with coordinates independent and identically distributed (iid, hereinafter) standard Normals, as in Muirhead (1982), but work with the covariance operator (instead of matrix) and characteristic function. This coordinate-free approach allows us to seamlessly subsume the possibility that the multivariate Normal distribution is supported on a proper subspace of $\Re^n$ which is not spanned by a subset of the standard orthonormal basis $\{e_1, \cdots, e_n\}$. In the spirit of Axler (2015), references to which in the sequel are by omission only, this relegates the manipulative matrix algebra that so dominates the multivariate Normal literature to the back burner.

The paper is organized as follows. In Section 2, we introduce all the notations, definitions, and results from linear algebra used in Section 3. The main result of this section is Theorem 1, where we show, using the Spectral Theorem, that the orthogonal projection operator onto the range of the adjoint of a linear operator $T$ is $UT$, where $U$ is an invertible linear operator. Theorem 1 is used in our proof of Theorem 3.

In section 3, we first summarize the basic facts about the multivariate Normal distribution in our coordinate-free setup. Theorem 2 shows that if $Y$ is $\Re^n$-valued multivariate Normal, then $(\mathcal{S}Y, \mathcal{T}Y)$ is $\Re^m \times \Re^p$-valued multivariate Normal, where $\mathcal{S} \in \mathcal{L}(\Re^n, \Re^m)$ and $\mathcal{T} \in \mathcal{L}(\Re^n, \Re^p)$. Corollary 1, which is used in our proof of Theorem 3, then formulates a necessary and sufficient condition in terms of the operators $D$, $\mathcal{S}$, and $\mathcal{T}$ for the independence of $\mathcal{S}Y$ and $\mathcal{T}Y$, where $D$ is the covariance of $Y$.

Note that, since a component of a vector is a linear transformation of the vector and a linear transformation of a multivariate Normal random variable is another multivariate Normal random variable [Lemma 5], Theorem 4 allows us to deduce Theorem 1.2.11(b) of Muirhead (1982) (and Proposition 3.13 of Eaton (1983)) as an immediate corollary.

We also present alternative derivations of the independence of the sample mean and the sample variance of a random sample from a Normal distribution [Remark 1], the "partial out" formula for population regression in the Normal model [Corollary 2], and the sufficiency of the sample mean in the Normal location model [Remark 3]. We simplify the expressions for the conditional mean and covariance obtained in Theorem 4 in Remark 2, leading to a direct verification of the iterated expectation and the analysis of variance formulae. We outline a direction in which our method can possibly be extended in Remark 5.

The Appendix contains two technical lemmas that are used in the proofs of Theorem 4 and Corollary 2.

The following notational conventions and their consequences are used throughout the rest of the paper. The equality of two random variables, unless otherwise mentioned, implies equality almost surely. For any Polish space $\mathfrak{X}$, let $\mathcal{B}(\mathfrak{X})$ denote the Borel $\sigma$-algebra of $\mathfrak{X}$.

Let $h$ be a map from an arbitrary set $\mathfrak{Y}$ into a measurable space $(\mathfrak{Z}, \mathcal{Z})$. For an arbitrary subset $B$ of $\mathfrak{Z}$, let $h^{-1}(B)$ denote $\{y \in \mathfrak{Y} : h(y) \in B\}$, whereas for an arbitrary subset $A$ of $\mathfrak{Y}$, let $h(A)$ denote $\{h(y) : y \in A\}$. Note that $h^{-1}(h(A)) = A$ for every subset $A$ of $\mathfrak{Y}$. Let $\sigma(h)$ denote the smallest $\sigma$-algebra of subsets of $\mathfrak{Y}$ that makes $h$ measurable. Since $\{h^{-1}(B) : B \in \mathcal{Z}\}$ is a $\sigma$-algebra of subsets of $\mathfrak{Y}$ [Dudley (1989, page 98)], we obtain $\sigma(h) = \{h^{-1}(B) : B \in \mathcal{Z}\}$.

**2. The results from Linear Algebra.** Let $V, W$ be finite-dimensional vector spaces. For any $\mathcal{S} \in \mathcal{L}(W, V)$, let $\mathcal{R}(\mathcal{S}) \subseteq V$ denote the range of $\mathcal{S}$ and $\mathcal{N}(\mathcal{S}) \subseteq W$ denote the null space of $\mathcal{S}$.

The result of Lemma 1 is used in the proof of Theorem 1. It is mentioned in Exercise 3.D.3; we present a proof for the sake of completeness.

**Lemma 1**. Let $V$ be a finite-dimensional vector space, $W$ a subspace of $V$, and $\mathcal{U} \in \mathcal{L}(W, V)$. There exists an invertible operator $U \in \mathcal{L}(V)$ such that $Uw = \mathcal{U}w$ for every $w \in W$ if and only if $\mathcal{U}$ is injective.

Proof of Lemma 1. Let $U \in \mathcal{L}(V)$ be invertible such that $Uw = \mathcal{U}w$ for every $w \in W$. Let $w \in W$ be such that $\mathcal{U}w = 0$, implying $Uw = 0$. Since $U$ is invertible, hence injective by Theorem 3.69, we conclude that $w = 0$, showing that $\mathcal{U}$ is injective.

To show the converse, let $Q$ be a direct sum complement of $W$ and $X$ a direct sum complement of $\mathcal{R}(\mathcal{U})$, where the existence of $Q$ and $X$ are guaranteed by Theorem 2.34. Let $\{q_1, \cdots, q_k\}$ be a basis for $Q$ and $\{x_1, \cdots x_m\}$ a basis for $X$. By the fundamental theorem of linear maps [Theorem 3.22], $\dim \mathcal{R}(\mathcal{U}) \leq \dim W$, implying, by Theorem 2.43, $k \leq m$. Recall that every $v \in V$ can be uniquely decomposed as $w + q$, where $w \in W$ and $q \in Q$. Define $U \in \mathcal{L}(V)$ as

$$Uv = \mathcal{U}w + x, \text{ where } x = \sum_{j=1}^{k} c_j x_j \in X \text{ and } q = \sum_{j=1}^{k} c_j q_j.$$

Since $X$ is a direct sum complement of $\mathcal{R}(\mathcal{U})$, $Uv = 0$ implies $\mathcal{U}w = x = 0$. Since $\mathcal{U}$ is injective and $\{x_1, \cdots x_k\}$ is linearly independent, $U$ is injective, hence invertible (again by Theorem 3.69). ☐

In the sequel we assume that $W$ and $V$ are real inner product spaces. For $\mathcal{T} \in \mathcal{L}(V, W)$, let $\mathcal{T}^* \in \mathcal{L}(W, V)$ denote the adjoint of $\mathcal{T}$ [Definition 7.2]. For any subspace $Q$ of $V$, let $Q^{\perp}$ denote the orthogonal complement of $Q$ [Definition 6.45] and $\Pi_Q \in \mathcal{L}(V)$ denote the orthogonal projection operator onto $Q$ [Definition 6.53].

**Theorem 1**. Given $T \in \mathcal{L}(V)$, there exists $U \in \mathcal{L}(V)$, invertible and depending on $T$, such that $\Pi_{\mathcal{R}(T^*)} = UT$.

Proof of Theorem 1. We first observe that $T^*T$ is a positive operator [Definition 7.31]. By Theorems 7.6(e) and 7.6(c), $T^*T$ is self-adjoint; the observation follows, since, by the definition of the adjoint operator, for every $v \in V$,

$$\langle T^*Tv, v \rangle = \|Tv\|^2 \geq 0. \tag{1}$$

By the Real Spectral Theorem [Theorem 7.29(b)], $V$ has an orthonormal basis $\{f_1, \cdots, f_n\}$ consisting of eigenvectors of $T^*T$ with corresponding eigenvalues $\{\lambda_1, \cdots, \lambda_n\}$. By Theorem 7.35(b), $\lambda_j \geq 0$ for all $1 \leq j \leq n$. Since, by $(1)$, $\lambda_j = \|Tf_j\|^2$, we obtain

$$\begin{aligned} \lambda_j = 0 \Leftrightarrow Tf_j = 0 \\ \lambda_j > 0 \Leftrightarrow Tf_j \neq 0. \end{aligned} \tag{2}$$

If $\lambda_j = 0$ for every $j = 1, \cdots, n$, then $T^*T$ is the zero operator, implying, by $(1)$, that $T$ is the zero operator. By Theorem 7.7, $T^*$ is the zero operator as well, and the theorem trivially holds with $U = I$. Hence, $\mathfrak{P} = \{j : 1 \leq j \leq n,\ \lambda_j > 0\}$ is non-empty without loss of generality; let $\mathfrak{P}^c = \{j : 1 \leq j \leq n,\ \lambda_j = 0\}$.

For $j \in \mathfrak{P}$, since $T^*\left(\frac{Tf_j}{\lambda_j}\right) = f_j$, we have $f_j \in \mathcal{R}(T^*)$, equivalently,

$$\Pi_{\mathcal{R}(T^*)} f_j = f_j. \tag{3}$$

For $j \in \mathfrak{P}^c$, $f_j \in \mathcal{N}(T)$ by $(2)$; since $\mathcal{N}(T) = (\mathcal{R}(T^*))^\perp$ [Theorem 7.7(c)],

$$\Pi_{\mathcal{R}(T^*)} f_j = 0. \tag{4}$$

By $(3)$ and $(4)$, for any $x \in V$,

$$\Pi_{\mathcal{R}(T^*)} x = \sum_{j \in \mathfrak{P}} \langle x, f_j \rangle f_j. \tag{5}$$

By $(2)$ and the definition of $\mathfrak{P}^c$, for any $x \in V$,

$$Tx = \sum_{j \in \mathfrak{P}} \langle x, f_j \rangle Tf_j. \tag{6}$$

By definition of $f_j$ and $\lambda_j$, $1 \leq j \leq n$, the list of vectors $\{Tf_j / \|Tf_j\| : j \in \mathfrak{P}\}$ in $V$ is orthonormal, and consequently, by Theorem 6.26, linearly independent; the same conclusion holds for the list $\{f_j : j \in \mathfrak{P}\}$. For $W = \text{span}\{Tf_j : j \in \mathfrak{P}\}$, $\mathcal{U} \in \mathcal{L}(W, V)$ defined by $\mathcal{U}Tf_j = f_j$ is clearly injective. By Lemma 1 there exists an invertible operator $U \in \mathcal{L}(V)$ such that

$$Uw = \mathcal{U}w \text{ for every } w \in W. \tag{7}$$

Now note that, for any $x \in V$, by $(6)$, $(7)$, the definition of $\mathcal{U}$, and $(5)$, in that order,

$$UTx = \sum_{j \in \mathfrak{P}} \langle x, f_j \rangle UTf_j = \sum_{j \in \mathfrak{P}} \langle x, f_j \rangle \mathcal{U}Tf_j = \sum_{j \in \mathfrak{P}} \langle x, f_j \rangle f_j = \Pi_{\mathcal{R}(T^*)} x,$$

completing the proof. $\square$

**Lemma 2**. Given $D \in \mathcal{L}(V)$ positive, there exists $D^{-1/2} \in \mathcal{L}(V)$ positive such that

$$D^{1/2} D^{-1/2} = D^{-1/2} D^{1/2} = \Pi_{\mathcal{R}(D)} \tag{8}$$

and

$$DD^{-1/2} = D^{1/2} = D^{-1/2} D, \tag{9}$$

where $D^{1/2}$ denotes the unique positive square root of $D$.

Proof of Lemma 2. By Theorem 7.36, $D^{1/2}$ exists and is defined by

$$D^{1/2} f_j = \lambda_j^{1/2} f_j, \tag{10}$$

where $\{f_1, f_2, \cdots, f_n\}$ is an orthonormal basis of $V$ consisting of eigenvectors of $D$ with corresponding eigenvalues $\{\lambda_1, \lambda_2, \cdots, \lambda_n\}$, and $\mathfrak{P}$ and $\mathfrak{P}^c$ are as in the proof of Theorem 1. Let $D^{-1/2} \in \mathcal{L}(V)$ be defined by

$$D^{-1/2} f_j = \begin{cases} \lambda_j^{-1/2} f_j & \text{if } j \in \mathfrak{P} \\ 0 & \text{if } j \in \mathfrak{P}^c. \end{cases} \tag{11}$$

Clearly, $\left\langle D^{-1/2} x, y \right\rangle = \sum_{j \in \mathfrak{P}} \lambda_j^{-1/2} \langle x, f_j \rangle \langle y, f_j \rangle = \left\langle x, D^{-1/2} y \right\rangle$, showing that $D^{-1/2}$ is self-adjoint and implying $\left\langle D^{-1/2} x, x \right\rangle \geq 0$, that is, $D^{-1/2} \in \mathcal{L}(V)$ is positive.

For any $y \in V$, by (10) and (11),

$$D^{1/2} D^{-1/2} y = \sum_{j \in \mathfrak{P}} \langle y, f_j \rangle f_j = D^{-1/2} D^{1/2} y. \tag{12}$$

Since $D$ is self-adjoint, by Theorem 7.7(d),

$$\mathcal{R}(D) = (\mathcal{N}(D))^{\perp}. \tag{13}$$

Since $\{f_j : j \in \mathfrak{P}^c\}$ is contained in $\mathcal{N}(D)$ and $\{f_j : j \in \mathfrak{P}\}$ in $\mathcal{R}(D)$, it follows that $\{f_j : j \in \mathfrak{P}\}$ is a basis of $\mathcal{R}(D)$, and (8) follows from (12). Note that (9) follows from the observation that both $DD^{-1/2} y$ and $D^{-1/2} D y$ equal $\sum_{j \in \mathfrak{P}} \langle y, f_j \rangle \lambda_j^{1/2} f_j = D^{1/2} y$. $\square$

**3. The multivariate Normal distribution results.** Let $Z_1, \cdots, Z_n$ be iid standard Normal random variables. The distribution of

$$Z = \sum_{k=1}^{n} e_k Z_k \tag{14}$$

is defined to be the standard multivariate Normal distribution $\mathfrak{N}_n(0, I)$, where $I \in \mathcal{L}(\Re^n)$ is the identity operator.

**Lemma 3**. The characteristic function $\Psi_{0,I}$ of the $\mathfrak{N}_n(0, I)$ distribution is given by

$$\Psi_{0,I}(t) = E(\exp(i\langle t, Z\rangle)) = \exp\Big( - \|t\|^2/2\Big).$$

Proof of Lemma 3. Follows from Proposition 9.4.2(a) of Dudley (1989). □

**Lemma 4**. Given $Z \sim \mathfrak{N}_n(0, I)$, $\mu \in \Re^n$, and $T \in \mathcal{L}(\Re^n)$, let

$$Y = \mu + TZ. \tag{15}$$

Then, for any $s, t \in \Re^n$,

$$E(\exp(i\langle t, Y\rangle)) = \exp\left( i\langle t, \mu\rangle - \frac{1}{2}\langle t, TT^* t\rangle \right)$$
$$E(\langle t, Y\rangle) = \langle t, \mu\rangle$$
$$\text{Cov}(\langle t, Y\rangle, \langle s, Y\rangle) = \big\langle t, TT^* s\big\rangle.$$

Proof of Lemma 4. Straightforward algebra using the bilinearity of the inner product and the definitions of $T^*$, $Z$ in (14), and $Y$ in (15), along with the fact that $Z_1, \cdots, Z_n$ are iid standard Normal random variables, proves the lemma. □

**Definition 1**. The distribution of $Y$ in (15) is defined to be the multivariate Normal distribution with mean $\mu$ and covariance $TT^*$. Recall that a distribution on $\Re^n$ is uniquely determined by its characteristic function [Dudley (1989, Theorem 9.5.1)]; since the characteristic function of $Y$ is determined by its mean $\mu$ and covariance $TT^*$, the multivariate Normal distribution $\mathfrak{N}_n(\mu, D)$ with mean $\mu$ and covariance $D$ (where $D$ is any positive operator on $\Re^n$) is uniquely defined in terms of the characteristic function

$$\Psi_{\mu,D}(t) = \exp\left( i\langle t, \mu\rangle - \frac{1}{2}\langle t, Dt\rangle \right). \tag{16}$$

**Lemma 5**. Given $Y \sim \mathfrak{N}_n(\mu, D)$ and $\mathcal{S} \in \mathcal{L}(\Re^n, \Re^m)$,

$$\mathcal{S}Y \sim \mathfrak{N}_m(\mathcal{S}\mu, \mathcal{S}D\mathcal{S}^*).$$

Proof of Lemma 5. The proof is a straightforward consequence of Definition 1. □

Given $(s_1, t_1), (s_2, t_2) \in \Re^m \times \Re^p = \Re^{m+p}$, the inner product on $\Re^m \times \Re^p$ is given by

$$\big\langle (s_1, t_1), (s_2, t_2)\big\rangle = \langle s_1, s_2\rangle + \langle t_1, t_2\rangle.$$

**Theorem 2**. Given $Y \sim \mathfrak{N}_n(\mu, D)$, $\mathcal{S} \in \mathcal{L}(\Re^n, \Re^m)$, and $\mathcal{T} \in \mathcal{L}(\Re^n, \Re^p)$, the random vector $(\mathcal{S}Y, \mathcal{T}Y) \sim \mathfrak{N}_{m+p}$ with mean $= (\mathcal{S}\mu, \mathcal{T}\mu) \in \Re^m \times \Re^p$ and covariance operator $\mathcal{K} \in \mathcal{L}(\Re^m \times \Re^p)$ given by

$$\mathcal{K}(s,t) = (\mathcal{S}D\mathcal{S}^*s + \mathcal{S}D\mathcal{T}^*t, \mathcal{T}D\mathcal{S}^*s + \mathcal{T}D\mathcal{T}^*t).$$

Proof of Theorem 2. We first verify that $\mathcal{K} : \Re^m \times \Re^p \mapsto \Re^m \times \Re^p$ is a positive operator. The verification of $\mathcal{K} \in \mathcal{L}(\Re^m \times \Re^p)$ being linear and self-adjoint is routine. Since

$$\begin{aligned}&\langle \mathcal{S}D\mathcal{S}^*s, s\rangle + \langle \mathcal{S}D\mathcal{T}^*t, s\rangle + \langle \mathcal{T}D\mathcal{S}^*s, t\rangle + \langle \mathcal{T}D\mathcal{T}^*t, t\rangle \\ &= \left\|D^{1/2}\mathcal{S}^*s\right\|^2 + 2\langle D^{1/2}\mathcal{T}^*t, D^{1/2}\mathcal{S}^*s\rangle + \left\|D^{1/2}\mathcal{T}^*t\right\|^2,\end{aligned}$$

the non-negativity of $\langle \mathcal{K}(s,t), (s,t)\rangle$ for every $(s,t) \in \Re^m \times \Re^p$ follows. By the definition of the inner product in $\Re^m \times \Re^p$ and (16), the characteristic function $\Psi$ of the random vector $(\mathcal{S}Y, \mathcal{T}Y)$ taking values in $\Re^m \times \Re^p$ is given by

$$\Psi(s,t) = \exp\left(i\langle \mathcal{S}^*s + \mathcal{T}^*t, \mu\rangle - \frac{1}{2}\langle \mathcal{S}^*s + \mathcal{T}^*t, D(\mathcal{S}^*s + \mathcal{T}^*t)\rangle\right). \tag{17}$$

Since $\langle \mathcal{S}^*s + \mathcal{T}^*t, \mu\rangle = \langle (s,t), (\mathcal{S}\mu, \mathcal{T}\mu)\rangle$ and

$$\begin{aligned}&\langle \mathcal{S}^*s + \mathcal{T}^*t, D(\mathcal{S}^*s + \mathcal{T}^*t)\rangle \\ &= \langle s, \mathcal{S}D\mathcal{S}^*s + \mathcal{S}D\mathcal{T}^*t\rangle + \langle t, \mathcal{T}D\mathcal{S}^*s + \mathcal{T}D\mathcal{T}^*t\rangle \\ &= \langle (s,t), (\mathcal{S}D\mathcal{S}^*s + \mathcal{S}D\mathcal{T}^*t, \mathcal{T}D\mathcal{S}^*s + \mathcal{T}D\mathcal{T}^*t)\rangle,\end{aligned} \tag{18}$$

the proof follows. □

**Corollary 1.** Given $Y \sim \mathfrak{N}_n(\mu, D)$, $\mathcal{S} \in \mathcal{L}(\Re^n, \Re^m)$, and $\mathcal{T} \in \mathcal{L}(\Re^n, \Re^p)$, $\mathcal{S}Y$ and $\mathcal{T}Y$ are independent if and only if $\mathcal{S}D\mathcal{T}^* \in \mathcal{L}(\Re^p, \Re^m)$, equivalently $\mathcal{T}D\mathcal{S}^* \in \mathcal{L}(\Re^m, \Re^p)$, is the zero operator.

Proof of Corollary 1. From (17),

$$\Psi(s,t) = \exp(i\langle s, \mathcal{S}\mu\rangle)\exp(i\langle t, \mathcal{T}\mu\rangle)\exp\left(-\frac{1}{2}\langle \mathcal{S}^*s + \mathcal{T}^*t, D(\mathcal{S}^*s + \mathcal{T}^*t)\rangle\right).$$

Now by (18),

$$\langle \mathcal{S}^*s + \mathcal{T}^*t, D(\mathcal{S}^*s + \mathcal{T}^*t)\rangle = \langle s, \mathcal{S}D\mathcal{S}^*s\rangle + 2\langle s, \mathcal{S}D\mathcal{T}^*t\rangle + \langle t, \mathcal{T}D\mathcal{T}^*t\rangle.$$

Thus, by (16) and Lemma 5,

$$\Psi(s,t) = E(\exp(i\langle s, \mathcal{S}Y\rangle))E(\exp(i\langle t, \mathcal{T}Y\rangle))\exp(-\langle s, \mathcal{S}D\mathcal{T}^*t\rangle).$$

That is, the characteristic function of $(\mathcal{S}Y, \mathcal{T}Y)$ is the product of two factors; one factor is the characteristic function of the product measure of the distributions induced by $\mathcal{S}Y$ on $\mathcal{B}(\Re^m)$ and $\mathcal{T}Y$ on $\mathcal{B}(\Re^p)$, whereas the other factor is $\exp(-\langle s, \mathcal{S}D\mathcal{T}^*t\rangle)$. Therefore, $\mathcal{S}Y$ and $\mathcal{T}Y$ are independent if and only if $\exp(-\langle s, \mathcal{S}D\mathcal{T}^*t\rangle) = 1$ for every $s \in \Re^m$ and $t \in \Re^p$, equivalently, $\langle s, \mathcal{S}D\mathcal{T}^* t \rangle = 0$ for every $s$ and $t$. Since $\mathcal{T}D\mathcal{S}^* = (\mathcal{S}D\mathcal{T}^*)^*$ by Theorems 7.6(e) and 7.6(c), the proof follows. $\square$

**Remark 1.** For $X_1, \cdots, X_n$ iid Normal random variables with mean $\theta$ and variance $\sigma^2$,

$$\bar{X} = \frac{1}{n}\sum_{i=1}^{n} X_i \text{ and } S^2 = \frac{1}{n-1}\sum_{i=1}^{n}(X_i - \bar{X})^2$$

are independent [Casella and Berger (2002, Theorem 5.3.1(a))]. Most textbooks prove this result by working with the (joint) density of the sample and using the Jacobian formula for finding the density of the transformation that maps the sample to the sample mean and the sample variance. Some textbooks use Basu's (1955) Theorem on an ancillary statistic (the sample variance) being independent of a complete sufficient statistic (the sample mean) to prove this result. We are going to show that this result is a straightforward consequence of Corollary 1.

Let $J$ be the sum of the standard orthonormal basis vectors, and $\{J\}$ the span of $J$. Since

$$\Pi_{\{J\}} x = \|J\|^{-2} \langle x, J \rangle J,$$

we have

$$\bar{X} = \|J\|^{-2} \langle \Pi_{\{J\}} X, J \rangle \text{ and } S^2 = \left(\|J\|^2 - 1\right)^{-1} \left\|\left(I - \Pi_{\{J\}}\right) X\right\|^2,$$

where $X = \sum_{k=1}^{n} X_k e_k \sim \mathfrak{N}_n(\theta J, \sigma^2 I)$. By Corollary 1, using the fact that an orthogonal projection operator is self-adjoint, we obtain $\Pi_{\{J\}} Y$ and $\left(I - \Pi_{\{J\}}\right) Y$ are independent if and only if $\left(I - \Pi_{\{J\}}\right) D \Pi_{\{J\}}$ is the zero operator, where $Y \sim \mathfrak{N}_n(\mu, D)$. Clearly, $\left(I - \Pi_{\{J\}}\right) D \Pi_{\{J\}} x = 0$ for all $x \in \Re^n$ if and only if $\left(I - \Pi_{\{J\}}\right) D J = 0$, equivalently, $J$ is an eigenvector of $D$. Since $J$ is an eigenvector of $\sigma^2 I$, the independence of $\bar{X}$ and $S^2$ follows from that of $\Pi_{\{J\}} X$ and $\left(I - \Pi_{\{J\}}\right) X$.

The class of positive operators with $J$ as an eigenvector does not reduce to the singleton set $\{I\}$. For $n = 2$, define $D \in \mathcal{L}(\Re^2)$ by $D(u_1, u_2) = (2u_1 + u_2, u_1 + 2u_2)$; clearly, $J = (1, 1)$ is an eigenvector for $D$. Thus, for the sample mean and the sample variance to be independent, it is not necessary for the sample to be iid. As long as the joint distribution of the sample is multivariate Normal such that $J$ is an eigenvector of the covariance, the independence of the sample mean and the sample variance holds. However, for Normal random variables that are dependent, whether the joint distribution is multivariate Normal becomes a modeling question, as the joint distribution of even pairwise uncorrelated Normal random variables may not be multivariate Normal. //

**Theorem 3**. Given $Y \sim \mathfrak{N}_n(\mu, D)$ and $T \in \mathcal{L}(\Re^n)$, define

$$S = TD^{1/2}; \tag{19}$$

then

$$D^{1/2}\Pi_{\mathcal{N}(S)}D^{-1/2}Y \text{ is independent of } TY \tag{20}$$

and

$$Y - D^{1/2}\Pi_{\mathcal{N}(S)}D^{-1/2}Y \text{ is an affine function of } TY. \tag{21}$$

Proof of Theorem 3. For any $x \in \Re^n$ and $F \in \mathcal{L}(\Re^n)$, we obtain from (9),

$$TD\big(D^{1/2}\Pi_{\mathcal{N}(F)}D^{-1/2}\big)^* x = TDD^{-1/2}\Pi_{\mathcal{N}(F)}D^{1/2}x = TD^{1/2}\Pi_{\mathcal{N}(F)}D^{1/2}x,$$

implying, by (19), that $TD\big(D^{1/2}\Pi_{\mathcal{N}(S)}D^{-1/2}\big)^*$ is the zero operator, whence (20) follows from Corollary 1.

To prove (21) we first observe that, by (13) and (8) in that order,

$$Y = \Pi_{\mathcal{R}(D)}Y + \Pi_{\mathcal{N}(D)}Y = D^{1/2}D^{-1/2}Y + \Pi_{\mathcal{N}(D)}Y. \tag{22}$$

Now we are going to show that

$$\Pi_{\mathcal{N}(D)}Y \;=\; \Pi_{\mathcal{N}(D)}\mu, \tag{23}$$

implying, by (22), that

$$Y = \; D^{1/2}D^{-1/2}Y + \Pi_{\mathcal{N}(D)}\mu. \tag{24}$$

Let $f_j$, $\lambda_j$, $\mathfrak{P}$, and $\mathfrak{P}^c$ be as in the proof of Lemma 2. Recall that $\{f_j : j \in \mathfrak{P}^c\}$ is an orthonormal basis for $\mathcal{N}(D)$ and $\{f_j : j \in \mathfrak{P}\}$ is an orthonormal basis for $\mathcal{R}(D)$. For $j \in \mathfrak{P}^c$ and $u \in \Re$, by (16),

$$E\big(\exp\big(iu\langle Y, f_j\rangle\big)\big) = \exp\bigg(i\langle uf_j, \mu\rangle - \frac{1}{2}\langle uf_j, Duf_j\rangle\bigg) = \exp\big(iu\langle \mu, f_j\rangle\big),$$

implying $\langle Y, f_j\rangle = \langle \mu, f_j\rangle$ for all $j \in \mathfrak{P}^c$, that is, (23).

By Theorems 6.47 and 7.7(c),

$$I = \Pi_{\mathcal{R}(S^*)} + \Pi_{\mathcal{N}(S)}, \tag{25}$$

implying, by (24),

$$Y - D^{1/2}\Pi_{\mathcal{N}(S)}D^{-1/2}Y \;= D^{1/2}\Pi_{\mathcal{R}(S^*)}D^{-1/2}Y + \Pi_{\mathcal{N}(D)}\mu. \tag{26}$$

By Theorem 1, there exists $U \in \mathcal{L}(\Re^n)$, invertible, such that

$$\Pi_{\mathcal{R}(S^*)} = US; \tag{27}$$

applying (27), (19), and (24) in that order, we obtain that RHS(26) equals

$$D^{1/2}USD^{-1/2}Y + \Pi_{\mathcal{N}(D)}\mu = D^{1/2}UTY + \big(I - D^{1/2}UT\big)\Pi_{\mathcal{N}(D)}\mu,$$

thereby establishing (21). □

**Theorem 4**. Given $Y \sim \mathfrak{N}_n(\mu, D)$ and $\mathcal{T} \in \mathcal{L}(\Re^n,\Re^m)$, the conditional distribution of $Y$ given $\mathcal{T}Y$ is multivariate Normal on $\mathcal{B}(\Re^n)$.

Proof of Theorem 4. The idea is to first construct $T \in \mathcal{L}(\Re^n)$ such that $\sigma(\mathcal{T}Y) = \sigma(TY)$ and then use Theorem 3 to find the conditional distribution of $Y$ given $TY$.

Let $\mathfrak{T} \in \mathcal{L}(\Re^n,\mathcal{R}(\mathcal{T}))$ be defined by $\mathfrak{T}x = \mathcal{T}x$ for every $x \in \Re^n$. Clearly, $\mathfrak{T}$ is surjective, hence injective and invertible [Theorem 3.69]. By Theorem 7.7, $\mathfrak{T}^* \in \mathcal{L}(\mathcal{R}(\mathcal{T}),\Re^n)$ is invertible as well. Let $T \in \mathcal{L}(\Re^n)$ be defined by $Tx = \mathfrak{T}^*\mathfrak{T}x$; note that $T$, being the composition of two invertible transformations, is invertible.

Let $(\Omega, \mathcal{F}, P)$ denote the probability space underlying $Y$. Since $\mathcal{T}$ is continuous, hence measurable, $\mathcal{T}Y : (\Omega, \mathcal{F}) \mapsto (\Re^m, \mathcal{B}(\Re^m))$ is measurable, implying $(\mathcal{T}Y)^{-1}(E) \in \mathcal{F}$ for every $E \in \mathcal{B}(\Re^m)$. Note that, $(\mathcal{T}Y)^{-1}(E) = (\mathfrak{T}Y)^{-1}(E \cap \mathcal{R}(\mathcal{T}))$. Further, $\mathcal{R}(\mathcal{T})$, by virtue of being closed, is an element of $\mathcal{B}(\Re^m)$, and $\mathcal{B}(\mathcal{R}(\mathcal{T}))$ is the trace of $\mathcal{B}(\Re^m)$ on $\mathcal{R}(\mathcal{T})$; that is, $\mathcal{B}(\mathcal{R}(\mathcal{T})) = \{E \cap \mathcal{R}(\mathcal{T}) : E \in \mathcal{B}(\Re^m)\}$, implying $\sigma(\mathcal{T}Y) = \sigma(\mathfrak{T}Y)$. Since $\mathfrak{T}$ and $T$ are injective (and measurable), by Lemma A.1, $\sigma(\mathfrak{T}) = \mathcal{B}(\Re^n) = \sigma(T)$; that is, $\{\mathfrak{T}^{-1}(G) : G \in \mathcal{B}(\mathcal{R}(\mathcal{T}))\} = \mathcal{B}(\Re^n) = \{T^{-1}(F) : F \in \mathcal{B}(\Re^n)\}$. Since, for any $G \in \mathcal{B}(\mathcal{R}(\mathcal{T}))$ and $F \in \mathcal{B}(\Re^n)$,

$$(\mathfrak{T}Y)^{-1}(G) = Y^{-1}\big(\mathfrak{T}^{-1}(G)\big) \text{ and } (TY)^{-1}(F) = Y^{-1}\big(T^{-1}(F)\big),$$

we obtain

$$\sigma(\mathfrak{T}Y) = \big\{(\mathfrak{T}Y)^{-1}(G) : G \in \mathcal{B}(\mathcal{R}(\mathcal{T}))\big\} = \big\{(TY)^{-1}(F) : F \in \mathcal{B}(\Re^n)\big\} = \sigma(TY).$$

That completes the proof that $\sigma(\mathcal{T}Y) = \sigma(TY)$.

By the pull-out property of conditional expectation, (20), and (21),

$$\begin{aligned} &E\Big(\exp\big(i\langle t, Y\rangle\big)\Big|TY\Big) \\ &= \exp\big(i\langle t, Y - D^{1/2}\Pi_{\mathcal{N}(S)}D^{-1/2}Y\rangle\big)E\big(\exp\big(i\langle t, D^{1/2}\Pi_{\mathcal{N}(S)}D^{-1/2}Y\rangle\big)\big), \end{aligned}$$

where $S$ is as in (19). By Lemma 5 and (16),

$$E\big(\exp\big(i\langle t, D^{1/2}\Pi_{\mathcal{N}(S)}D^{-1/2}Y\rangle\big)\big) = \exp\bigg(i\langle t, D^{1/2}\Pi_{\mathcal{N}(S)}D^{-1/2}\mu\rangle - \frac{1}{2}\langle t, Gt\rangle\bigg),$$

where the operator $G = D^{1/2}\Pi_{\mathcal{N}(S)}D^{-1/2}DD^{-1/2}\Pi_{\mathcal{N}(S)}D^{1/2}$ equals, by (9) and (8), $D^{1/2}\Pi_{\mathcal{N}(S)}\Pi_{\mathcal{R}(D)}\Pi_{\mathcal{N}(S)}D^{1/2}$. Therefore, the conditional distribution of $Y$ given $TY$, hence that given $\mathcal{T}Y$, is Normal with mean $\upsilon(Y) = Y - D^{1/2}\Pi_{\mathcal{N}(S)}D^{-1/2}(Y - \mu)$ and covariance $G$. □

**Remark 2.** The expressions for the mean $\upsilon$ and the covariance $G$ of the conditional distribution of $Y$ given $\mathcal{T}Y$ can be considerably simplified, rendering the verification of the iterated expectation formula and the analysis of variance formula immediate.

Since $Y - \mu = \Pi_{\mathcal{R}(D)}(Y - \mu)$ by $(23)$, applying $(25)$ and $(8)$ in that order,

$$\upsilon(Y) = \mu + D^{1/2}\Pi_{\mathcal{R}(S^*)}D^{-1/2}(Y - \mu). \tag{28}$$

To simplify the expression of the covariance operator $G$ (which does not depend on $Y$), first note that, by $(8)$ and $(9)$,

$$D^{1/2}\Pi_{\mathcal{R}(D)} = D^{1/2}D^{1/2}D^{-1/2} = DD^{-1/2} = D^{1/2}. \tag{29}$$

By $(27)$, $(19)$, and $(29)$, in that order,

$$\Pi_{\mathcal{R}(S^*)}\Pi_{\mathcal{R}(D)} = UTD^{1/2}\Pi_{\mathcal{R}(D)} = UTD^{1/2} = \Pi_{\mathcal{R}(S^*)}. \tag{30}$$

Consequently, $G$ equals

$$\begin{aligned} & D^{1/2}\Pi_{\mathcal{R}(D)}\Pi_{\mathcal{N}(S)}D^{1/2} - D^{1/2}\Pi_{\mathcal{R}(S^*)}\Pi_{\mathcal{R}(D)}\Pi_{\mathcal{N}(S)}D^{1/2} && \text{by } (25) \\ &= D^{1/2}\Pi_{\mathcal{N}(S)}D^{1/2} - D^{1/2}\Pi_{\mathcal{R}(S^*)}\Pi_{\mathcal{R}(D)}\Pi_{\mathcal{N}(S)}D^{1/2} && \text{by } (29) \\ &= D^{1/2}\Pi_{\mathcal{N}(S)}D^{1/2} - D^{1/2}\Pi_{\mathcal{R}(S^*)}\Pi_{\mathcal{N}(S)}D^{1/2} && \text{by } (30) \\ &= D^{1/2}\Pi_{\mathcal{N}(S)}D^{1/2} && \text{by } (25). \end{aligned} \tag{31}$$

Recall that if $V$ and $W$ are random vectors in $\Re^n$ such that $V$ is an affine function of $W$, i.e., $V = a + RW$, where $a \in \Re^n$ and $R \in \mathcal{L}(\Re^n)$, then

$$E(V) = a + R(E(W)) \text{ and } \mathrm{Cov}(V) = R\mathrm{Cov}(W)R^*,$$

implying $E\big(E\big(Y\big|TY\big)\big) = \mu$, verifying the iterated expectation formula.

Defining the expected value of an operator-valued random element $\Upsilon$ as the operator $E(\Upsilon)$ such that $E(\Upsilon)(s) = E(\Upsilon(s))$ for every $s \in \Re^n$, provided the expected value on the right hand side is well defined for every $s$, we conclude from $(31)$ that

$$E\big(\mathrm{Cov}\big(Y\big|TY\big)\big) = D^{1/2}\Pi_{\mathcal{N}(S)}D^{1/2}.$$

By $(28)$ and Lemma 5,

$$\mathrm{Cov}\big(E\big(Y\big|TY\big)\big) = D^{1/2}\Pi_{\mathcal{R}(S^*)}D^{-1/2}DD^{-1/2}\Pi_{\mathcal{R}(S^*)}D^{1/2} = D^{1/2}\Pi_{\mathcal{R}(S^*)}D^{1/2},$$

where the second equality follows by $(9)$, $(8)$, and $(30)$. The analysis of variance formula, $E\big(\mathrm{Cov}\big(Y\big|TY\big)\big) + \mathrm{Cov}\big(E\big(Y\big|TY\big)\big) = \mathrm{Cov}(Y)$, is verified by $(25)$. //

An immediate corollary of Theorem 4 is the "partial out" formula for population regression in the Normal model.

**Corollary 2.** If $(X, Y, Z)$ is multivariate Normal such that $Z = (Z_1, \cdots, Z_{n-2}) \in \Re^{n-2}$ and $X, Y \in \Re$, then

$$E(Y|X,Z) - E(Y|Z) = \frac{\text{Cov}(X,Y|Z)}{\text{Var}(X|Z)}[\text{Var}(X|Z) > 0](X - E(X|Z)). \quad (32)$$

Proof of Corollary 2. Let $\mathcal{P}_{1,3} \in \mathcal{L}(\Re^n, \Re^{n-1})$ and $\mathcal{P}_3 \in \mathcal{L}(\Re^n, \Re^{n-2})$ be defined by $\mathcal{P}_{1,3}w = (x, z)$ and $\mathcal{P}_3 w = z$, where $w = (x, y, z)$. By (28),

$$\text{LHS}(32) = \left\langle D^{1/2}\left(\Pi_{\mathcal{R}(S_{1,3}^*)} - \Pi_{\mathcal{R}(S_3^*)}\right)D^{-1/2}(W - \mu), e_2\right\rangle, \quad (33)$$

where $W = (X, Y, Z)$, $\mu \in \Re^n$ is the mean of $W$ and $D \in \mathcal{L}(\Re^n)$ is the covariance of $W$, $S_{1,3} = P_{1,3}D^{1/2}$ with $P_{1,3} \in \mathcal{L}(\Re^n)$ defined by $P_{1,3}w = (x, 0, z)$, and $S_3 = P_3 D^{1/2}$ with $P_3 \in \mathcal{L}(\Re^n)$ defined by $P_3 w = (0, 0, z)$.

To show $\text{RHS}(32) = \text{RHS}(33)$, we first observe, using (31) and (28),

$$\begin{aligned} \text{Cov}(X,Y|Z) &= \langle D^{1/2}\Pi_{\mathcal{N}(S_3)}D^{1/2}e_1, e_2\rangle \\ \text{Var}(X|Z) &= \langle D^{1/2}\Pi_{\mathcal{N}(S_3)}D^{1/2}e_1, e_1\rangle = \left\|\Pi_{\mathcal{N}(S_3)}D^{1/2}e_1\right\|^2 \\ X - E(X|Z) &= \left\langle \left(I - D^{1/2}\Pi_{\mathcal{R}(S_3^*)}D^{-1/2}\right)(W - \mu), e_1\right\rangle. \end{aligned} \quad (34)$$

Since $P_{1,3}e_j = e_j$ if $j \neq 2$, $P_{1,3}e_2 = 0$, $P_3 e_j = e_j$ if $j \geq 3$, and $P_3 e_j = 0$ if $j = 1, 2$,

$$\begin{aligned} \mathcal{N}(S_{1,3}) &= \left(\text{span}\left\{D^{1/2}e_1, D^{1/2}e_3, \cdots, D^{1/2}e_n\right\}\right)^\perp \\ \mathcal{N}(S_3) &= \left(\text{span}\left\{D^{1/2}e_3, \cdots, D^{1/2}e_n\right\}\right)^\perp. \end{aligned}$$

If $D^{1/2}e_1 \in \text{span}\left\{D^{1/2}e_3, \cdots, D^{1/2}e_n\right\} = \mathcal{R}(S_3^*)$, then $\text{RHS}(33) = 0$, in which case $\text{Var}(X|Z)$, a constant function of $Z$, is also equal to $0$, thereby vacuously satisfying (32).

If $D^{1/2}e_1 \notin \text{span}\left\{D^{1/2}e_3, \cdots, D^{1/2}e_n\right\}$, by Lemma A.2,

$$\begin{aligned} &\left(\Pi_{\mathcal{R}(S_{1,3}^*)} - \Pi_{\mathcal{R}(S_3^*)}\right)D^{-1/2}(W - \mu) \\ &= \left\|\Pi_{\mathcal{N}(S_3)}D^{1/2}e_1\right\|^{-2}\langle \Pi_{\mathcal{N}(S_3)}D^{-1/2}(W - \mu), \Pi_{\mathcal{N}(S_3)}D^{1/2}e_1\rangle \Pi_{\mathcal{N}(S_3)}D^{1/2}e_1. \end{aligned}$$

Therefore, RHS (33) equals

$$\frac{\langle \Pi_{\mathcal{N}(S_3)}D^{-1/2}(W - \mu), \Pi_{\mathcal{N}(S_3)}D^{1/2}e_1\rangle}{\left\|\Pi_{\mathcal{N}(S_3)}D^{1/2}e_1\right\|^2}\langle D^{1/2}\Pi_{\mathcal{N}(S_3)}D^{1/2}e_1, e_2\rangle,$$

which, by the first two rows in (34), equals, since $\Pi_{\mathcal{N}(S_3)}$ is self-adjoint and idempotent,

$$\frac{\text{Cov}(X,Y|Z)}{\text{Var}(X|Z)}\big\langle D^{1/2}\Pi_{\mathcal{N}(S_3)}D^{-1/2}(W-\mu),e_1\big\rangle.$$

Since $D^{1/2}\Pi_{\mathcal{N}(S_3)}D^{-1/2}=I-D^{1/2}\Pi_{\mathcal{R}(S_3^*)}D^{-1/2}-\Pi_{\mathcal{N}(D)}$ by (25), (8), and (13), the proof follows by the third row in (34) and (23). □

**Remark 3.** Given a random sample $X_1,\cdots,X_n$ from the Normal distribution with mean $\theta$ and (known) variance $\sigma^2$, $\bar{X}$ is a sufficient statistic for $\theta$ [Example 6.2.4, Casella and Berger (2002)]. While the typical proof uses the powerful factorization theorem [Theorem 6.2.6, Casella and Berger (2002)], we are going to show that the conditional distribution of the sample $X=\sum_{k=1}^{n}X_ke_k\sim\mathfrak{N}_n(\theta J,\sigma^2 I)$ given $\bar{X}$, that is, $n^{-1}\Pi_{\{J\}}X$, does not depend on $\theta$, thereby proving the sufficiency of $\bar{X}$ directly. The said conditional distribution, by Theorem 2, (28), and (31), is multivariate Normal with mean $\theta J+D^{1/2}\Pi_{\mathcal{R}(S^*)}D^{-1/2}(X-\theta J)$ and covariance $D^{1/2}\Pi_{\mathcal{N}(S)}D^{1/2}$, where $D=\sigma^2 I$ and $S=n^{-1}\Pi_{\{J\}}D^{1/2}$, implying $S^*=\sigma n^{-1}\Pi_{\{J\}}$ and consequently, $\Pi_{\mathcal{R}(S^*)}=\Pi_{\{J\}}$, further implying that the mean equals $\Pi_{\{J\}}X$ and the covariance equals $\sigma^2\big(I-\Pi_{\{J\}}\big)$. //

**Remark 4.** In a companion paper we are currently working on, we use the *local linearity* of a $C^1$ vector field $g$ (on $\Re^n$ that takes values in $\Re^m$) and the decomposition of $Y$ in Theorem 3 to approximate the conditional distribution of $Y$ given $g(Y)$. If $g$ is constant, then the $\sigma$-algebra generated by $g(Y)$ is the trivial $\sigma$-algebra consisting of the empty set and the entire sample space, making $Y$ independent of $g(Y)$, so that the conditional distribution of $Y$ given $g(Y)$ is the unconditional distribution of $Y$. If $g$ is injective, $\sigma(g)$ equals $\mathcal{B}(\Re^n)$ by Lemma A.1; consequently, the conditional distribution of $Y$ given $g(Y)$ is the point mass at $Y$.

Thus, the interesting problem unfolds when $g$ is a non-constant, non-injective, $C^1$ vector field. The conditional distribution of $Y$ given $g(Y)$ will not be multivariate Normal if $g$ is not linear. Let $\mathfrak{g}:\Re^n\mapsto\mathcal{L}(\Re^n,\Re^m)$ be defined to be the continuous function that maps $a\in\Re^n$ to the total derivative of $g$ at $a$. Given $\epsilon>0$, there exists a compact subset $K_\epsilon$ of $\Re^n$ such that $P(Y\in K_\epsilon)>1-\epsilon$. The function $\mathfrak{g}$, restricted to $K_\epsilon$, is uniformly continuous, and admits a uniformly continuous modulus of continuity $\zeta_\epsilon$. We are working to show that, for a suitably chosen metric for weak convergence of probability measures and $\delta>0$ (that depends on the given $\epsilon$ through the supremum of the modulus of continuity $\zeta_\epsilon$ on the closed interval $[0,B_\epsilon]$, where $B_\epsilon$ is the diameter of $K_\epsilon$), the conditional distribution of $Y$ given $g(Y)$ is within $\delta$-neighborhood of the family of multivariate Normal distributions. //

**Remark 5.** Proposition 3.13 of Eaton (1983) holds at a much greater level of generality; see Bogachev (1998, Theorem 3.10.1). Extending our approach, in the absence of an inner product, to finding the conditional distribution of a Gaussian random element given a (non-injective, non-trivial) linear transformation appears to be an interesting area of future research. In particular, if $Y$ is the Brownian motion, i.e., the random element in $\mathcal{C}([0,1])$ distributed according to the Wiener measure, $\nu_1, \cdots, \nu_k$ are measures on $\mathcal{B}([0,1])$, and $\mathcal{T} : \mathcal{C}([0,1]) \mapsto \Re^k$ is given by $\mathcal{T}f = (\int f d\nu_1, \cdots, \int f d\nu_k)$, what is the conditional distribution of $Y$ given $\mathcal{T}Y$? //

**Appendix.** We present a lemma on measurability [Lemma A.1] that was used in the proof of Theorem 4 and another lemma on orthogonal projection operators in a Hilbert space [Lemma A.2] that was used in the proof of Corollary 2.

**Lemma A.1** Let $\mathfrak{X}$ and $\mathfrak{Y}$ be Polish spaces. If $h : \mathfrak{X} \mapsto \mathfrak{Y}$ is Borel measurable and injective, then $\sigma(h) = \mathcal{B}(\mathfrak{X})$.

Proof of Lemma A.1. The inclusion $\sigma(h) \subseteq \mathcal{B}(\mathfrak{X})$ follows from the pertinent definitions. For $M \in \mathcal{B}(\mathfrak{X})$, define $B = h(M)$. Since $h^{-1}(B) = M$ and $B \in \mathcal{B}(\mathfrak{Y})$ by Theorem I.3.9 of Parthasarathy (1967), the reverse inclusion follows. □

**Lemma A.2** For a proper and closed subspace $V$ of a Hilbert space $H$ and $x \in H \setminus V$, let $V_x$ denote the closure of the subspace $\{v + cx : v \in V, c \in \Re\}$. Then, for any $y \in H$,

$$\Pi_{V_x} y - \Pi_V y = \left\| \Pi_{V^\perp} x \right\|^{-2} \left\langle \Pi_{V^\perp} y, \Pi_{V^\perp} x \right\rangle \Pi_{V^\perp} x. \tag{35}$$

Proof of Lemma A.2. By the definition of $V_x$ there exists a sequence $\{v_n : n \geq 1\} \subset V$ and a sequence $\{c_n : n \geq 1\} \subset \Re$ such that

$$\Pi_{V_x} y = \lim_{n \to \infty} (v_n + c_n x). \tag{36}$$

Since $V \subset V_x$, $\Pi_V y = \Pi_V \Pi_{V_x} y$; since $\Pi_V$ is continuous, by (36),

$$\Pi_V y = \lim_{n \to \infty} (v_n + c_n \Pi_V x). \tag{37}$$

Since $x - \Pi_V x = \Pi_{V^\perp} x$, subtracting (37) from (36),

$$\text{LHS}(35) = \left( \lim_{n \to \infty} c_n \right) \Pi_{V^\perp} x. \tag{38}$$

Since $\text{LHS}(35) = \Pi_{V^\perp} y - \Pi_{(V_x)^\perp} y$, taking the inner product of both sides of (38) with $\Pi_{V^\perp} x$ and using the linearity and homogeneity of inner product, we obtain

$$\left\langle \Pi_{V^\perp} y, \Pi_{V^\perp} x \right\rangle - \left\langle \Pi_{(V_x)^\perp} y, \Pi_{V^\perp} x \right\rangle = \left( \lim_{n \to \infty} c_n \right) \left\| \Pi_{V^\perp} x \right\|^2. \tag{39}$$

Since $x \in V_x$ and $\Pi_V x \in V \subset V_x$, $\Pi_{V^\perp} x \in V_x$, implying $\left\langle \Pi_{(V_x)^\perp} y, \Pi_{V^\perp} x \right\rangle = 0$; since $x \notin V$, that is, $\left\| \Pi_{V^\perp} x \right\|^2 \neq 0$, the lemma follows from (39) and (38). □

RAJESHWARI MAJUMDAR
rajeshwari.majumdar@uconn.edu
PO Box 47
Coventry, CT 06238

SUMAN MAJUMDAR
suman.majumdar@uconn.edu
1 University Place
Stamford, CT 06901